\documentclass[twoside,11pt]{amsart}
\usepackage{amsmath}
\usepackage{amssymb}
\usepackage{latexsym}
\usepackage{verbatim}
\usepackage{graphics}
\usepackage{graphicx}

\newtheorem{theorem}{Theorem}
\newtheorem{corollary}{Corollary}
\newtheorem{lemma}{Lemma}
\newtheorem{proposition}{Proposition}

\newtheorem{definition}{Definition}
\newtheorem{remark}{Remark}

\def\Z{{\mathbb Z}}

\def\G{{\bf G}}
\def\S{{\mathbb S}}
\def\F{{\mathbb F}}
\def\K{{\mathbb K}}
\def\ZZ{{\mathbb Z}}
\def\P{{\mathcal P}}

\begin{document}

\title {$G$-odometers and their almost 1-1 extensions.}

\author{Mar\'{\i}a Isabel Cortez, Samuel Petite}
\address{Centro de Modelamiento Matem\'atico,  Av. Blanco Encalada
2120 Piso 7
Santiago, Chile. \\
Laboratoire d'Analyse, Topologie, Probabilit{\'e}s, CNRS U.M.R. 6632, Universit\'e Paul C\'ezanne, 13
397 Marseille Cedex 20  France.} \email{mcortez@dim.uchile.cl, samuel.petite@univ.u-3mrs.fr}

\subjclass{Primary: 54H20; Secondary: 37B50 } \keywords{Almost 1-1 extensions, odometer, Toeplitz system,
discrete group actions}

\maketitle{}
\begin{abstract}
In this paper we recall the concepts of $G$-odometer and $G$-subodometer for $G$-actions, where $G$ is a
discrete finitely generated group, which generalize the notion of odometer in the case $G=\ZZ$. We
characterize the $G$-regularly recurrent systems as the minimal almost 1-1 extensions of subodometers, from
which we deduce that the family of the $G$-Toeplitz subshifts coincides with the family of the minimal
symbolic almost 1-1 extensions of subodometers.
\end{abstract}
\section{Introduction}

It is well known that any continuous dynamical system factorizes onto a minimal equicontinuous dynamical
system (see \cite{A}). For this reason, it is useful to study the minimal equicontinuous factors of a general
dynamical system, or conversely, to determine the extensions of a particular minimal equicontinuous system.
We place us in this latest problematic. The aim of this paper is to characterize the extensions of a
particular type of equicontinuous factor: the $G$-odometers, where $G$ is a  discrete finitely generated
group, like for example a non Abelian free group. The notion of $G$-odometer
generalizes the notion of odometer in the case $G=\Z$.\\
An example of extensions of odometers are the Toeplitz flows, which
were introduced by Jacobs and Keane in \cite{JK}. Toeplitz flows
have been extensively studied in different contexts and they have
been used to provide a series of examples with interesting dynamical
properties (see for example \cite{D}, \cite{GJ}, \cite{Wi}).
 Markley and Paul characterize them in \cite{MP} as the minimal almost 1-1 extensions of
odometers and  a proof of this theorem is given in \cite{DL} by Downarowicz and Lacroix.\\
Following the work developed in \cite{Co} for $G = \Z^d$, we prove   that for a discrete finitely generated
group $G$, the $G$-Toeplitz systems are the symbolic minimal almost 1-1 extensions of $G$-odometers. The main
difference that appears with the Abelian case, is the existence of some
degenerated systems that we call {subodometers}.\\
This paper is organized as follows: in Section 2, we give some basic definitions relevant for the
study of topological dynamical systems. We recall also the generalized notions of odometer and
subodometer and we identify the set of eigenvalues of these systems. In Section 3, we introduce the
notions of regularly recurrent systems and  strongly regularly recurrent systems. We characterize
them as the minimal almost 1-1
extensions of subodometers and  odometers respectively. 
In the particular case where $G$ is amenable, we show in Section 4
that the set of invariant probability measures of a $G$-regularly
recurrent Cantor system can be represented as an inverse limit. In
Section 5, in the case when $G$ is a residually finite group, we
introduce a notion of semicocycles and we show that an almost 1-1
extension of a $G$-subodometer is conjugated to the action of $G$ on
some semicocycle. Finally in Section 6, we consider a particular
family for a discrete group $G$: the $G$-Toeplitz arrays, which is a
particular family of semicocycles when $G$ is residually finite. We
prove, by giving an explicit construction, that this family
coincides with the family of symbolic almost 1-1 extensions of the
$G$-subodometers.


\medskip

{\bf Acknowledgment: }it is a great pleasure for S.P. to thank the Instituto de Matem\'aticas
U.N.A.M. of Cuernavaca (M{\'e}xico), where part of this work has been done, for their warm hospitality.

\section{Basic definitions and background}
In this article,  by a  {\em topological dynamical system} we mean a pair  $(X, G)$, where $G$ is a
topological group which acts, by homeomorphism, on a compact metric space $(X,d)$. Given $g\in G$
and $x\in X$ we will identify $g$ with the associated homeomorphism and we denote by $g.x$ the
action of $g$ on $x$. The dynamical system $(X, G)$ is {\em free} if $g.x=x$ for some $x\in X$
implies $g=e$, where $e$ is the neutral element in $G$. For a syndetic subgroup $\Gamma$ of $G$,
the {\em $\Gamma$-orbit} of $x\in X$ is $O_{\Gamma}(x)=\{\gamma.x: \gamma\in \Gamma\}$ and the {\em
$\Gamma$-system associated} to $x$  is $(\Omega_{\Gamma}(x),\Gamma)$, where $\Omega_{\Gamma}(x)$ is
the closure of $O_{\Gamma}(x)$ and the action of $\Gamma$ on $\Omega_{\Gamma}(x)$ is the
restriction to $\Gamma$ and $\Omega_{\Gamma}(x)$ of the action of $G$ on $X$. The set of {\em
return times} of $x\in X$ to $A\subseteq X$ is $ T_{A}(x)=\{g\in G:\, g.x\in A\}$.  The topological
dynamical system $(X, G)$ is {\em minimal} if the orbit of any $x\in X$  is dense in $X$, and it is
said to be {\em equicontinuous} if for every $\varepsilon>0$ there exists $\delta>0$ such that if
$x,y\in X$ satisfy $d(x,y)<\delta$ then $d(g.x, g.y)<\varepsilon$ for all $g\in G$. We say that
$(X, G)$ is an {\em extension} of $(Y, G)$, or that $(Y, G)$ is a {\em factor} of $(X,G)$, if there
exists  a continuous surjection $\pi:X\to Y$ such that $\pi$ preserves the  action. We call $\pi$ a
{\em factor map}. When the factor map is bijective, we say that $(X, G)$ and $(Y, G)$ are {\em
conjugate}. The factor map $\pi$ is an {\em almost 1-1 factor map} and $(X, G)$ is an {\em almost
1-1 extension} of $(Y, G)$ by $\pi$ if the set of points having one pre-image  is residual
(contains a dense $G_{\delta}$ set) in $Y$. In the minimal case it is equivalent to the existence
of a point with one pre-image.\\
The set $\mathcal M_G(X)$ of {\em invariant probability measures} of
$X$ is the set of probability measures $\mu$ defined on $\mathcal
B(X)$, the Borel $\sigma$-algebra of $X$, such that
$\mu(g.B)=\mu(B)$ for all $g\in G$ and $B\in \mathcal{B}(X)$.

\subsection{$G$-odometers and $G$-subodometers}
In all the following, we will denote by $G$ a discrete group generated by a finite family and by
$e$ its neutral element.
\begin{definition}\label{defRF}
A discrete finitely generated group $G$ is called residually finite if and only if there exists a
sequence $\Gamma_1 \supset \Gamma_2 \supset \ldots \supset \Gamma_n \supset \ldots$ of subgroups
$\Gamma_n$ with finite index in $G$ such that:
$$\bigcap_n \Gamma_n = \{e\}.$$
\end{definition}

A trivial example of a residually finite subgroup is the group of integers $\Z$, for example by
taking the groups $\Gamma_n = n! \Z$. Less trivial examples are given by the fundamental groups of
connected oriented compact graph. When $\pi: S_2 \to S_1$ is a finite covering of an oriented compact connected
graph $S_2$ onto a compact graph $S_1$, the application $\pi$ induces an homomorphism $\pi_*$
from the fundamental group of $S_2$ to the fundamental group of $S_1$. The image of the morphism
$\pi_*$ is a subgroup of the fundamental group of $S_1$. The index of this subgroup is then the
number of pre-image of one point for the map $\pi$. Let us denote by $\widetilde{S_1}$ the
universal cover of $S_1$. Consider a sequence $(S_n, \pi_n)_n$ of finite covering $\pi_n :
S_{n+1}\to S_n$ of compact conneced and oriented graph $S_n$ such that for each $n$ the injectivity radius of
$\widetilde{S_1}$ onto $S_n$ goes to infinity when $n$ goes to infinity. The sequence of
fundamental groups of graphs $S_n$ satisfies then the condition of the Definition \ref{defRF}. More
generally, we have the following result of Mal'cev \cite{Ma}:
\begin{theorem}{\cite{Ma}}
For any integer $n$ and any  field $\K$ with characteristic null, every subgroup finitely generated
of the group of invertible  matrices $GL(n, \K)$ is a residually finite groups.
\end{theorem}

\noindent In particular, the free groups $\F_n$ with $n$ generators, the groups of surfaces and the
braids group $B_n$ generated by $n$ elements are residually finite groups.

Let us denote, for a subgroup $H$ of $G$, by $G/H$ the set of right class of $H$ in $G$. It is important to note that $G$ acts on $G/H$ by left multiplication on the $H$-class. Now we will prove the useful following lemma:
\begin{lemma}\label{lemfond}
Let $G$ be a group. If $H$ is a subgroup of $G$ with index in $G$
equal to $n$ (i.e. the cardinal of the quotient space $G/ H$ is $n$)
then there exists a normal subgroup $K$ of $G$ included in $H$ such
that $G/K$ divide $n!$.
\end{lemma}
\begin{proof}
The group $G$ acts on $G/ H$ by left multiplication. This action defines an homomorphism $\rho$
from $G$ to the permutation group of $n$ elements. The kernel of this application is a normal
subgroup of $G$ included in $H$ and its index in $G$ divides the cardinal of permutations of
$n$ elements.
\end{proof}

\noindent As a corollary, $G$ is a residually finite group if and only if there exists a sequence
$H_1 \supset \ldots \supset H_n \supset \ldots$ of normal subgroups of $G$ with finite index in $G$
such that $\bigcap_n H_n =\{e\}$. Remark that, up to consider a quotient space, all finitely
generated groups are residually finites.

\bigskip

Let us consider a  discrete group $G$ generated by a finite family  and a decreasing sequence (for the
inclusion) $(\Gamma_i)_{i\geq 0}\subseteq G$ of subgroups with finite index in $G$ (we do not ask
$\bigcap_{i\geq 0}\Gamma_i=\{e\}$)  and let $\pi_i:G/\Gamma_{i+1}\to G/\Gamma_i$ be the function induced by
the inclusion $\Gamma_{i+1}\subset \Gamma_i$, $i\geq 0$. Consider the inverse limit
$$\overleftarrow{G}=\lim_{\leftarrow i}(G/\Gamma_i,\pi_i).$$

More precisely, $\overleftarrow{G}$ is defined as the subset of the product $\Pi_{i\geq
0}G/\Gamma_i$ consisting of the elements $\textbf{g}=(g_i)_{i\geq 0}$ such that
$\pi_i(g_{i+1})=g_i$ for all $i\geq 0$.

Every $G/\Gamma_i$ is endowed with the discrete topology and $\Pi_{i\geq 0}G/\Gamma_i$ with the
product topology. Thus $\overleftarrow{G}$ is a compact metrizable space whose topology is spanned
by the cylinder sets
$$[i;a]=\{\textbf{g}\in \overleftarrow{G}: g_i=a \},\mbox{ with $a \in
G/\Gamma_i$ and $i\geq 0$}.$$

\noindent The space $\overleftarrow{G}$ is a totally disconnected, it is a Cantor set when $G/ \cap_{i\geq 0}
\Gamma_i$ is infinite and a finite set when $G/ \cap_{i\geq 0} \Gamma_i$ is finite.

The group $G$ acts continuously on $\overleftarrow{G}$ by left multiplication,
namely for ${\bf g}
= (g_i)_i \in \overleftarrow{G}$ and $h \in G$,
$$ h. {\bf g}= (h._i g_i)_i, $$
\noindent where $h._i$ denotes the  action on $G/ \Gamma_i$ given by $h._ig\Gamma_i=hg\Gamma_i$,
for every $h\in G$ and $g\in G$.

\noindent Since for all $h \in G$ and for all cylinders $[i;a]$ we have
$${h}.([i;a])\subseteq [i;h._i a_i],$$ the topological dynamical system
$(\overleftarrow{G},G)$ is equicontinuous. Moreover, every orbit for this action is dense, then
$(\overleftarrow{G},G)$ is a minimal equicontinuous system.

\begin{definition}
We call $(\overleftarrow{G},G)$ a {\it $G$-subodometer system}\footnote{Note that this definition is not a
profinite completion of the group $G$ because here, we consider only a sequence of decreasing subgroups.} or
simply a {\it subodometer}. If in addition, every $\Gamma_i$ is normal, we say that $(\overleftarrow{G},G)$
is a  {\it $G$-odometer system} or simply an {\it odometer}.
\end{definition}

\noindent It is straightforward to show that for a point ${\bf g} = (g_i)_i$ of a subodometer
$\overleftarrow{G}$, its stabilizer for the $G$-action is the group $\bigcap_i \tilde{g}_i \Gamma_i
\tilde{g}^{-1}_i$, where $\tilde{g}_i$ is a representing element of the class $g_i\in \G/\Gamma_i$ in $G$,
for $i\geq 0$. Hence, when $G$ is a residually finite group and $\bigcap_{i\geq 0}\Gamma_i=\{e\}$, for ${\bf
e}$ the element $(e_i)_i$ of the $G$-subodometer $\overleftarrow{G}$, where $e_i$ is the projection of the
neutral element of $G$ on $G/\Gamma_i$, the stabilizer of $\bf e$ is trivial. This does not mean necessarily
that the action of $G$ on $\overleftarrow{G}$ is free. If furthermore, all the groups $\Gamma_i$ are normal
subgroups of $G$, then the stabilizer of every point of a $G$-odometer is trivial and the action of $G$ is
free. For this reason we call, when $G$ is residually finite and $\bigcap_{i\geq 0}\Gamma_i=\{e\}$, the
$G$-odometer $\lim_{\leftarrow i} (G/\Gamma_i, \pi_i)$ a {\em free $G$-odometer}.

\bigskip

\noindent  If $(\overleftarrow{G},G)$ is an odometer then the set $\overleftarrow{G}$ is a group equipped
with the multiplication defined by
$$\textbf{g}.\textbf{h}=(g_i._i h_i)_{i\geq 0},  $$

\noindent  where $._i$ denotes the multiplication operation induced on $G/\Gamma_i$ by the multiplication on
$G$. Remark that for a free odometer $(\overleftarrow{G}, G)$, the group $G$ is then a dense subgroup of
$\overleftarrow{G}$.

\noindent Notice that for all ${\bf g}$ in a cylinder set $[i;a]$ of an odometer $\overleftarrow{G}
=\lim_{\leftarrow i}(G/H_i,\pi_i)$, the set of return times of $\bf g$ to $[i;a]$ is $H_i$. Through
this paper we will use this property and we will identify $\overleftarrow{G}$ with
$(\overleftarrow{G}, G)$.

 \begin{lemma}\label{factor}
Let $\overleftarrow{G}_j=\lim_{\leftarrow i}(G/H_i^{j},\pi_i)$ be two subodometers ($j=1,2$). Let
$\bf e_j$ $(j=1, 2)$ be the element  $(e^j_i)_i \in \overleftarrow{G}_j$ where
$e^j_i$ denotes the class of the neutral element $e \in G$ in $G/H_i^{j}$.\\
There is a factor map $\pi: (\overleftarrow{G}_1,G)\to (\overleftarrow{G}_2,G)$ such that $\pi
({\bf e_1})= {\bf e_2}$ if and only if for every $H_i^{2}$ there exists some $H_k^{1}$ such that
$H_k^{1}\subseteq H_i^{2}$.
\end{lemma}

\begin{proof}
If $\pi:\overleftarrow{G}_1\to \overleftarrow{G}_2$ is a factor map then by continuity, given
$i\geq 0$ and $e^2_i$ in $G/H_i^{2}$, there exists $k\geq 0$  such that $[k;e^1_k]\subseteq
\pi^{-1}[i;e^2_i]$. Let $v\in H_k^{1}$, we have that $v.\textbf{g}\in [k;e^1_k]$ for all
$\textbf{g}\in [k;e^1_k]$, which implies that
$$\pi(v.\textbf{g})=v.\pi(\textbf{g})\in [i;e^2_i].$$ Since
$\pi(\textbf{g})\in [i;e^2_i]$ and $T_{[i;e^2_i]}(\pi(\textbf{g}))=H_i^{2}$, we
get $v\in H_i^{2}$.\\
Suppose that for every $i\geq 0$ there exists $H_{n_i}^{1}\subseteq H_i^{2}$. Since the sequences
$(H_{i}^{j})_{i\geq 0}$, $j=1,2$, are decreasing, we can take $n_{i}\leq n_{i+1}$ for all $i\geq
0$. The function  $\pi: \overleftarrow{G}_1\to \overleftarrow{G}_2$ defined by $\pi((g_i)_{i\geq
0})=(j_{n_i}(g_{n_i}))_{i\geq 0}$ where $j_{n_i}: G/H_{n_i}^{1}\to G/H_i^{2}$ is the function
induced by the inclusion $H_{n_i}^{1}\subseteq H_i^{2}$, is a factor map.
\end{proof}

\medskip

By a straightforward application of the former lemma and Lemma \ref{lemfond}, we get
\begin{proposition}\label{propodo}
If $(\lim_{\leftarrow i}(G/\Gamma_i,\pi_i), G)$ is a $G$-subodometer,
then there
exists a
$G$-odometer which is an extension of this subodometer.
\end{proposition}



\begin{proposition}\label{caractersubodometre}

For $\overleftarrow{G}$ a $G$-odometer and $(X,G)$ a dynamical system, if there exists a factor
map from $\overleftarrow{G}$ onto $X$, then there exists a closed subgroup $H$ of $\overleftarrow{G}$ such that the
dynamical system $(\overleftarrow{G}/H, G)$ is conjugated to $(X,G)$.
\end{proposition}
\noindent In particular this proposition says that a subodometer is conjugate to the quotient of an odometer
by a closed subgroup.

\begin{proof}
Let us denote by $p$ the factor map $\overleftarrow{G} \to X$, $\bf e$ the neutral element of
$\overleftarrow{G}$ and for ${\bf g}$ an element of $\overleftarrow{G}$, we denote by $(g_i)_i$ a
sequence of $G \subset \overleftarrow{G}$ that converges to ${\bf g}$. Let $H$ be the closed subset
$p^{-1}(p(e))$ of $\overleftarrow{G}$. For ${\bf g}$ and ${\bf h}$ in $H$, we have:
$$p({\bf h}{\bf g})= \lim_i p( h_i g_i)= \lim_i h_i.p(g_i) =\lim_i
h_i.p(e)=\lim_i p(h_i)= p(e). $$ \noindent With the same technic we get:
$$p(({\bf g})^{-1})= \lim p(g_i^{-1} e)= \lim_i g_i^{-1}.p(e)=\lim_i
g_i^{-1}.p(g_i)=p(e).$$
\noindent  So ${\bf g}{\bf h}$ and ${\bf g}^{-1}$ belong to $H$, and
$H$ is a group.\\
Now let us see that $p^{-1}(p({\bf g}))= {\bf g}H$ for any ${\bf g} \in
\overleftarrow{G}$. Let
${\bf h}$ be in $H$, we have:
$$p({\bf g}{\bf h})= \lim_i p(g_i h_i) = \lim_i g_i.p(h_i)=\lim_i g_i.p(e)=
p({\bf g}).$$
\noindent Then ${\bf g}H \subset p^{-1}(p({\bf g}))$.\\
Let ${\bf h}\in \overleftarrow{G}$ be such that $p({\bf h})=p({\bf
g})$. Then $\lim_i p(g_i)=\lim_i p(h_i)$ and $p(e)=\lim_i
g_i^{-1}.h_i.p(e) = p({\bf g}^{-1}{\bf h})$. So ${\bf g}^{-1}{\bf
h}$ belongs to $H$ and $p^{-1}(p({\bf g}))={\bf g}H$. Therefore, the
map $p$ factorizes onto a homeomorphism from $\overleftarrow{G}/H$
to $X$.
\end{proof}

\subsection{Eigenvalues of odometers and subodometers.}
Let $(X,\mu,G)$ be a measure-theoretic dynamical system with a left
action of $G$. A character $\chi$ is a homomorphism from $G$ to the
group $\S^1$, the set of complex numbers with module 1.
Since the group $G$ is equipped with the discrete topology, every
character is a continuous map.\\
A character is an {\em eigenvalue} of $X$ if there exists $f\in
L^2_{\mu}(X)\setminus \{0\}$ such that $f(g.x)=\chi(g)f(x)$ for all
$x\in X$ and $g\in G$. We call $f$ an {\em eigenfunction} associated
to $\chi$. We say that an eigenvalue is a {\em continuous
eigenvalue} if it has an associated
continuous eigenfunction.\\


Since a $G$-odometer $\overleftarrow{G}$ is a compact group, the
normalized Haar measure left invariant $\lambda$ of
$\overleftarrow{G}$ is the only invariant probability measure of
$\overleftarrow{G}$ for the action of $G$. A $G$-subodometer is a
factor of a $G$-odometer (Proposition \ref{propodo}), therefore, by
this factorization, the subodometer inherits also of an invariant
measure for the $G$-action.  Since any factor map extends to a
continuous affine onto map between the set of invariant
probability measures \cite{DGS}, we conclude that a subodometer is
uniquely ergodic. \noindent Thus when we speak about a subodometer
$\overleftarrow{G}$ as a measure-theoretic dynamical system, we mean
$\overleftarrow{G}$ equipped with the only invariant probability
measure $\lambda$   for the action of $G$.
\begin{proposition}
\label{eigenvalue} Let $\overleftarrow{G}=\lim_{\leftarrow
n}(G/\Gamma_n,\pi_n)$ be a subodometer. The set of eigenvalues of
$\overleftarrow{G}$ is given by $ E_G=\bigcup_{n\geq 0}\{character \
\chi :G \to \S^1, \ \chi(\gamma)=1 \ {\textrm{for }\ \gamma \in
\Gamma_n}\}.$ Moreover, every eigenvalue of $\overleftarrow{G}$ is a
continuous eigenvalue.
\end{proposition}
\begin{proof}
For $n\geq 0$ we call $C_n=[n;e]$. Since $v,w\in G$ satisfy
$v.C_n=w.C_n$ if and only if  $w$ and $v$ belong to the same class
in $G/\Gamma_n$, it makes sense to write $v.C_n$ for $v\in
G/\Gamma_n$. Notice that the collection $\mathcal{P}_n=\{v.C_n :
v\in
G/\Gamma_n\}$ is a clopen partition of $G$.\\
Let $\chi \in E_G$ and let $n\geq 0$ be  such that $\chi(\gamma)=1$
for all $\gamma\in \Gamma_n$. This means that $\chi$ is constant on
each class of $G/\Gamma_n$, which implies that $f=\sum_{v\in
G/\Gamma_n}\chi(v)1_{v.C_n}$ is a well defined continuous function
that verifies $f(h.\textbf{g})=\chi(h)f(\textbf{g})$ for all
$\textbf{g}\in \overleftarrow{G}$ and
$h\in G$.\\
Let $\chi$ be an eigenvalue of $\overleftarrow{G}$ and let $f\in
L^2_{\lambda}(\overleftarrow{G})\setminus \{0\}$ be an associated
eigenfunction.
For $g\in G$ we have that
$$
\chi(g)\left(\int_{C_n}fd\lambda\right)=\int_{g.C_n}fd\lambda.
$$
Since $C_n=\gamma.C_n$ for all $\gamma\in \Gamma_n$, it holds that
\begin{equation}
\label{int}
\chi(g)\left(\int_{C_n}fd\lambda\right)=\int_{C_n}fd\lambda\, \, \,
\, \mbox{ for all $g\in \Gamma_n$}.
\end{equation}
Observe that
$$
\mathbb E(f|\mathcal{P}_n)=\sum_{g\in
K_n}\frac{\chi(g)}{\lambda(C_n)}\left(\int_{C_n}fd\lambda\right)
1_{g.C_n},
$$
for a finite set $K_n\subset G$ containing  at least one element of
each class of $G/\Gamma_n$. Since
$\mathcal{B}(\mathcal{P}_n)\uparrow \mathcal{B}(\overleftarrow{G})$,
by the increasing Martingale theorem, we have that $\mathbb
E(f|\mathcal{P}_n)$ converges to $f$ in
$L^2_{\lambda}(\overleftarrow{G})$. Because $f\neq 0$, this implies
there exists $m\geq 0$ such that $\int_{C_m}fd\lambda\neq 0$ and, by
(\ref{int}), we conclude that $\chi(\gamma)=1$ for all $\gamma\in
\Gamma_m$, which means that $\chi \in E_G$.
\end{proof}

\section{Characterization of minimal almost 1-1 extensions of odometers }

Let $(X,G)$ and $(Y,G)$ be two topological dynamical systems. $(Y,G)$ is said to be the {\em
maximal equicontinuos factor} of $(X,G)$ if it is an equicontinuos factor  of $(X,G)$ such that for
any other equicontinuous factor $(Y',G)$ of $(X,G)$ there exists a factor map $\pi:Y\to Y'$ that
satisfies $\pi
\circ f=f'$, with $f:X\to Y$ and $f':X\to Y'$  factor maps.\\
It is well known that every topological dynamical system has a maximal
equicontinuous factor and if
$(X,G)$ is a minimal almost 1-1 extension of a minimal equicontinuous system
$(Y,G)$, then $(Y,G)$
is the maximal equicontinuous factor of $(X,G)$ (for more details see
\cite{A}).

\subsection{Regularly recurrent systems.}

A subset $S$ of $G$ is said to be {\em syndetic} if there exists a compact subset $K$ of $G$ such
that $G=K.S=\{k.s:s\in S, \, k\in K\}$. Because we consider  $G$ a  discrete group, a subset $S$ of
$G$ is syndetic if and only there exists a finite subset $K$ of $G$ such that $G=K.S$. It is important to note that a subgroup $\Gamma$ of $G$ is syndetic if and only if $G/\Gamma$ is finite.\\
Let $(X,G)$ be a topological dynamical system and let $x\in X$. The point $x$ is {\em uniformly
recurrent} if for every open neighborhood $V$ of $x$ the set $T_V(x)$ is syndetic.
It is well known that $(\Omega_{G}(x),G)$ is
minimal if and only if $x$ is uniformly recurrent.\\
A point $x\in X$ is {\em regularly recurrent} if for every open
neighborhood $V$ of $x$ there is a syndetic subgroup $\Gamma$ of $G$
such that $\Gamma\subseteq T_V(x)$. We say that a
system is {\em regularly recurrent} if it is the orbit closure of a regularly
recurrent point.\\
Similarly, we say that a point $x\in X$ is {\em strongly regularly recurrent} if for every open
neighborhood $V$ of $x$ there is a clopen subset $W\subset V$, neighborhood of $x$, such that
$T_W(x)$ is  a  syndetic normal subgroup of $G$. We say that a system is {\em strongly regularly
recurrent} if it is the orbit closure of a strongly regularly recurrent point. Obviously, a
strongly regularly recurrent point is a regularly recurrent point. Regularly recurrent systems are
minimal.


The subodometers (resp. odometers) are examples of (resp. strongly) regularly recurrent systems. Moreover,
every point in a subodometer (resp. odometer) $\overleftarrow{G}$ is regularly recurrent (resp. strongly
regularly recurrent).

In this section, we will show that (resp. strongly) recurrent systems are exactly the minimal almost 1-1
extensions of the subodometers (resp. odometers). From that we will conclude that a group $G$ admits an
action that is strongly regularly recurrent and free if and only if $G$ is residually finite.

\begin{lemma}\label{minimal}
Let $(X,G)$ be a minimal topological dynamical system and let $x\in
X$. If $\Gamma\subseteq G$ is a syndetic subgroup  of $G$ then
$(\Omega_\Gamma(x),\Gamma)$ is minimal.
\end{lemma}
\begin{proof}
Let $H$ be a normal subgroup of $G$ included in $\Gamma$ (Lemma \ref{lemfond}). The group $G$ acts
by the natural product action on the compact spaces $X \times G/H$ and $X\times G/\Gamma$. Pick a
minimal set $M$ in $X\times G/H$. This set projects onto a minimal subset of $X$ hence onto $X$.
Thus for every $x\in X$ there exists a point $(x, a)\in M$ and this point is uniformly recurrent.
The right multiplication by $a^{-1}$ on the second axis is a conjugacy that sends the minimal set
$M$ onto a minimal set $M'$ that contains $(x,e)$. This set projects onto a minimal set of $X
\times G/\Gamma$ that contains the point $(x,[e])$ where $[e]$ denotes the $\Gamma$-class of the
neutral element $e$. This implies that for any neighborhood $V\subseteq X$ of $x$, the set $\{g:
g.x\in V, \, g\in \Gamma\}$ is syndetic.
\end{proof}

\begin{lemma}
\label{existencia} Let $(X,G)$ be a topological dynamical system and
let $x\in X$ be a (resp. strongly) regularly recurrent point. For
every closed neighborhood $V$ of $x$ there exists a (resp. normal)
syndetic subgroup $\Gamma$ of $G$   such that $\Gamma \subseteq
T_V(x)$ and $\{w(\Omega_\Gamma(x))\}_{w\in G/\Gamma}$ is a clopen
partition of $X$.
\end{lemma}
\begin{proof}
Let $V\subseteq X$ be a closed neighborhood of a  point $x$ regularly recurrent and let
$\Gamma'\subseteq G$ be a subgroup with finite index such that $\Gamma'\subseteq T_V(x)$. Let us
consider the normal subgroup $H \subset \Gamma'$ given by Lemma \ref{lemfond}. By Lemma
\ref{minimal}, the set $\Omega_H (x)$ is a closed set minimal and invariant for the $H$-action.
Since $H$ is normal, for any $g \in G$, the set $g. \Omega_H (x)$, which equals $\Omega_H(g.x)$, is
also closed, invariant and minimal for the $H$-action. Therefore if $w.\Omega_H(x)\cap
u.\Omega_H(x)\neq \emptyset$ for $u,w\in G$, we have  $w.\Omega_H(x)= u.\Omega_H(x)$.

\noindent  Furthermore, if $u$ and $w \in G$ are in the same $H$-class,
then  we
have also
$w.\Omega_H(x)= u.\Omega_H(x)$. Since $H$ is syndetic and the $G$-orbit of $x$
is dense, we have
$X=\bigsqcup_{u\in K} u.\Omega_H(x)$, for some finite family $K$ of $G$.

Let $\Gamma$ be the group $$\Gamma =\{ g \in G: \ g.\Omega_H(x) =\Omega_H (x)
\}.$$

\noindent We have $H \subset \Gamma$, so $\Gamma$ is syndetic. Since $\Omega_\Gamma (x)=
\Omega_H(x)$, we have $\Gamma \subset T_V(x)$ and  for any $g \in G$ $g.\Omega_\Gamma(x)$ and
$\Omega_\Gamma(x)$ are disjoint or equal  because they are minimal closed $H$-invariant sets. Thus
we get :
\begin{enumerate}
\item $g.\Omega_\Gamma (x)= g'.\Omega_\Gamma (x)$ {\textrm if and only if}  $
g\in g'\Gamma.$

\item $T_{g.\Omega_\gamma(x)}(y)=g\Gamma g^{-1}$ for every $y \in
g.\Omega_\gamma(x)$.
\end{enumerate}

\noindent It holds that for $w \in G/\Gamma$, $w.\Omega_\Gamma(x)$ is well
defined and
$\{w.\Omega_\Gamma(x)\}_{w \in G/\Gamma}$ is a clopen partition of $X$.

\medskip

When $x$ is a strongly regularly recurrent point of $X$,  we do exactly the same proof with $H$
being the normal subgroup $T_W(x)$ given by a  clopen neighborhood $W\subset V$ of $x$. Thanks this
strong property, we have that the group $\Gamma$ is actually the group $H$ and thus $\Gamma$ is a
normal subgroup of $G$.
\end{proof}

\begin{corollary}
\label{caracterizacion}

Let $(X,G)$ be a topological dynamical system and let $x\in X$. The
point $x$ is
(resp. strongly)
regularly recurrent if and only if there exists $(C_i)_{i\geq 0}$, a
fundamental system of clopen
neighborhoods of $x$ ($\cap_i C_i = \{x\}$), such that  for all $y\in C_i$ the
set of return times
of $y$ to $C_i$ is a syndetic (resp. normal) subgroup $\Gamma_i$ of $G$, for
every $i\geq 0$.
\end{corollary}

\begin{proof}
If  $x\in X$ has a fundamental system of neighborhoods as written
above, it is a (resp. strongly) regularly recurrent point.\\

\noindent The sequences $(C_i)_i$ and  $(\Gamma_i)_i$ are defined by induction. If $x$ is a
(resp. strongly) regularly recurrent point, let $C_1$ be the space $X$ and $\Gamma_1$ be the group
$G$.

\noindent So, given $C_n$ and $\Gamma_n$, we take an open neighborhood $V_{n+1}$ of $x$, whose the
closure is strictly included in  $C_n$. By Lemma \ref{existencia}, we obtain a syndetic (resp.
normal) group $\Gamma_{n+1}$ with $\Gamma_{n+1}\subseteq T_{\overline{V}_{n+1}}(x)$ and
$\{w(\Omega_{\Gamma_{n+1}}(x))\}_{w\in G/\Gamma_{n+1}}$ is a clopen partition of $X$. Clearly, we
have $\Gamma_{n+1}\subset \Gamma_n$. We set $C_{n+1}=\Omega_{\Gamma_{n+1}}(x)$ which is a clopen
set with $T_{C_{n+1}}(y)=\Gamma_{n+1}$ for all $y\in \Gamma_{n+1}$.

\noindent Since $\lim_{i\to \infty}{\rm diam}(V_n)=0$, we obtain that $(C_i)_{i\geq 0}$ is a
fundamental system of clopen neighborhoods of $x$.
\end{proof}

\bigskip

\begin{theorem}\label{almost}
A minimal topological dynamical system $(X,G)$ is an almost 1-1 extension of a subodometer (resp. odometer)
$\overleftarrow{G}$ by $\pi$ if and  only if $(X,G)$ is a  (resp. strongly) regularly recurrent system.
Moreover, the set of (resp. strongly) regularly recurrent points of $X$ is exactly the pre-image of the set
of points in $G$ which have only one pre-image by $\pi$.
\end{theorem}
\begin{proof}
Let $(X,G)$ be a minimal 1-1 extension of an  subodometer $\overleftarrow{G}=\lim_{\leftarrow
i}(G/\Gamma_i,\pi_i)$ (resp. odometer). Let $\pi: X\to \overleftarrow{G}$ be the almost 1-1 factor map and
let $x\in X$ be such that  $\{x\}= \pi^{-1} (\{\pi(x)\})$. Since $\pi$ is continuous,  if
$\pi(x)=(a_i)_{i\geq 0}\in \overleftarrow{G}$ then $(\pi^{-1}([i;a_i]))_{i}$ is a decreasing sequence of
clopen neighborhoods of $x$ that satisfies $$\bigcap_{i\geq 0}\pi^{-1}([i;a_i])=\{x\}.$$ We know that for
every $\textbf{g}\in [i;a_i]$, the set $T_{[i;a_i]}(\textbf{g})$ is a group conjugated to $\Gamma_i$,
therefore for all $y$ in $\pi^{-1}([i;a_i])$, we have $T_{\pi^{-1}([i;a_i])}(y)$ is a group conjugated to
$\Gamma_i$. So, by Corollary \ref{caracterizacion} we conclude that $x$ is a (resp. strongly)
regularly recurrent point of $X$.\\

\noindent Let $X$ be a (resp. strongly) regularly recurrent system and let $x\in X$ be a (resp. strongly)
regularly recurrent point with a trivial stabilizer. By Corollary \ref{caracterizacion}  there exists a
decreasing sequence $(C_i)_{i\geq 0}$ of clopen neighborhoods of $x$ such that $\bigcap_{i\geq 0}C_i=\{x\}$,
and there is a syndetic (resp. normal) subgroup $\Gamma_i$ such that $T_{C_i}(y)=\Gamma_i$ for all $y\in
C_i$, $i\geq 0$. Since $C_{i+1}\subseteq C_i$, we have that $\Gamma_{i+1}\subseteq \Gamma_i$, $i\geq 0$. So,
we can define the subodometer (resp. odometer) $\overleftarrow{G}=\lim_{\leftarrow i}(G/\Gamma_i,\pi_i)$. We
define $\pi:X\to \overleftarrow{G}$ by $\pi=(f_i)_{i\geq 0}$ where $f_i$ is the continuous map $f_i:X\to
G/\Gamma_i$ given by $f_i(y)=[z]$, where $[z]$ denotes the $\Gamma_i$-class of $z\in G$, if and only if $y\in
z.C_i$ for $y\in X$, $z\in G$ and $i\geq 0$. The function $\pi$ is a factor map, and, since $\bigcap_{i\geq
0}C_i=\{x\}$, we have that $\pi^{-1}\{\textbf{e}\}=\{x\}$. So, $\pi$ is an almost 1-1 extension.\\
If $\pi':X\to \overleftarrow{G'}$ is another almost 1-1 factor map and $\overleftarrow{G'}$ an
subodometer (resp. odometer), $\overleftarrow{G}$ and $\overleftarrow{G'}$ are the maximal
equicontinuous factor of $(X,G)$, therefore, they are conjugate. Thus there exists a factor map
$\pi'':\overleftarrow{G'}\to \overleftarrow{G}$ such that $\pi''\circ \pi'=\pi$, which implies that
$\pi'^{-1}\{x\}=\pi^{-1}\{\pi''(x)\}$ for any $x$ of $\overleftarrow{G'}$. We conclude that the set
of (resp. strongly) regularly recurrent points is exactly the pre-image of the points in $G$ which
have only one pre-image.
\end{proof}

By a straightforward application of Theorem \ref{almost} we get the following
corollaries.
\begin{corollary}
Every point of a (resp. strongly) regularly recurrent system $(X,G)$ is (resp. strongly) regularly recurrent
if and only if $(X,G)$ is conjugate to a (resp. odometer) subodometer.
\end{corollary}

\begin{corollary}
A discrete group finitely generated $G$ admits a strongly regularly recurrent free action on a compact metric
space if and only if $G$ is residually finite.
\end{corollary}

\begin{corollary}
Let $(X,G)$ be a regularly recurrent system and let
$\overleftarrow{G}$ be its maximal equicontinuous factor. The set of
continuous eigenvalues of $X$ is $E_G$.
\end{corollary}
\begin{proof}
It is clear that $E_G$ is contained in the set of continuous eigenvalues of $X$. Conversely, if
$\chi$ is a continuous eigenvalue of $X$ we can take $f:X\to \S^1$ an associated continuous
eigenfunction which is a factor map between $(X,G)$ and the dynamical system $(f(X),G)$, where the
action of $g\in G$ on $\exp(2i\pi x)\in f(X)$ is given by $g.\exp(2i\pi x)=\chi(g)\exp(2i\pi x)$,
which is an isometry. Thus the system $(f(X),G)$ is equicontinuous and therefore there exists a
factor map $\pi:\overleftarrow{G}\to f(X)$. Since $\pi$ is an eigenfunction associated to $\chi$ we
conclude that $\chi\in E_G$.
\end{proof}

\section{Regularly recurrent Cantor systems with $G$ amenable.}

We say that a topological dynamical system $(X,G)$ is a {\it (resp.
strongly) regularly recurrent Cantor system} if it is (resp.
strongly) regularly recurrent and $X$ is a Cantor set. In this
section we suppose that $(X,G)$ is a regularly recurrent Cantor
system.

\begin{proposition}
\label{KRpartition} Let $(X,G)$ be a regularly recurrent Cantor
system. There exists a sequence
$$
\mathcal{P}_n=\{w.C_{n,k}: w\in D_n, \, 1\leq k\leq k_n\},$$
  of finite clopen partitions of $X$ satisfying, for every $n\geq 0$, the
following:
\begin{enumerate}
\item
$C_{n+1}\subseteq C_n=\bigcup_{k=1}^{k_n}C_{n,k} \subset X$.
\item
There exists a syndetic subgroup $\Gamma_n$ of $G$ such that $D_n$
is a subset of $G$ containing exactly one representing element of
each class in $G/\Gamma_n$ and such that $T_{C_n}(x)=\Gamma_n$, for
all $x\in C_n$.
\item
$\mathcal{P}_{n+1}$ is finer than $\mathcal{\P}_n$.
\item
The collection of set $(P_n)_{n\geq 0}$ spans the topology of $X$.
\end{enumerate}

\end{proposition}

\begin{proof}
The idea of the proof is the same used in \cite{HPS} and \cite{Pu} to show that any
minimal Cantor $\ZZ$-system has a nested sequence of clopen Kakutani-Rohlin partitions.\\
We recall the algorithm introduced in \cite{Pu} to generate a
Kakutani-Rohlin partition finer that another one.
Let $\mathcal{R}$ be a finite clopen partition of $X$. Suppose that
$$ \mathcal{Q}=\{w.C_j: w\in D, \, 1\leq j\leq k\},$$
is another clopen partition of $X$  for which
 $k<\infty$,
 there exists a syndetic subgroup $\Gamma$ of $G$ such that
 $D=\{w_1,\cdots, w_l\}$ is a subset of $G$ containing exactly one
representing element of each class in $G/\Gamma$, and   the set of return times of any point in
$C=\bigcup_{j=1}^kC_j$ to $C$ is equal to $\Gamma$.  The next algorithm
produce a partition $\mathcal{R}\wedge\mathcal{Q}=\{w.B_j: w\in D,\, 1\leq j\leq d\}$  verifying
\begin{itemize}
  \item $\mathcal{R}\wedge\mathcal{Q}$ is finer than $\mathcal{R}$ and
$\mathcal{Q}$.
\item
$C=\bigcup_{j=1}^dB_j$
\end{itemize}

{\it \underline{Step 1:}} let $1\leq j\leq k$.  Consider
$A_{1,j,i_1},\cdots, A_{1,j,i_{l_{1,j}}}$,  the sets in $\mathcal{R}$ such
that
$$
w_1^{-1}.A_{1,j,i_s}\cap C_j\neq \emptyset,  \mbox{ for every }
1\leq s\leq l_{1,j}.$$
We denote by $B_{1,1},\cdots, B_{1,k_1}$, with  $k_1=\sum_{j=1}^{k}l_{1,j}$, the elements of the collection
$$\{w_1^{-1}.A_{1,j,i_s}\cap C_j: 1 \leq s\leq l_{1,j},1 \leq j\leq k\,\}.$$ We have that
$\mathcal{Q}_1=\{w.B_{1,j}: w\in D, \, 1\leq j\leq k_1\}$  is a clopen finite partition of $X$.
In addition,  for every $1\leq i\leq k_1$  there exist  $1\leq j\leq k$ and $1\leq s\leq l_{1,j}$ such that
 $w_1.B_{1,i}\subseteq A_{1,j,i_s}$, $w_1.B_{1,i}\subseteq
w_1.C_j$ and  $\bigcup_{s=1}^{l_j}B_{1,i}=C_j$. In other words, we have obtained a clopen
partition  $ \mathcal{Q}_1=\{w.B_{1,j}: w\in D, \, 1\leq j\leq k_1\},$
satisfying
\begin{itemize}
\item
For every $1\leq j\leq k_1$, there exist  $A$ in $\mathcal{R}$ and
$B$ in $\mathcal{Q}$ such that $w_1.B_{1,j}$ is contained in $A\cap B$.
\item
$\bigcup_{j=1}^{k_1}B_{1,j}=C.$
\end{itemize}
Now, for $2\leq n\leq l$, we suppose that the step $n-1$ has
produced a finite clopen partition $\mathcal{Q}_{n-1}=\{w.B_{n-1,j}: w\in
D,\, 1\leq j\leq k_{n-1}\}$ such that
\begin{itemize}
\item
For every $1\leq j \leq k_{n-1}$ and every $1\leq i\leq n-1$, there
exists $A$ in $\mathcal{R}$ and $B\in \mathcal{Q}$ such that
$w_i.B_{n-1,j}$ is contained in  $A\cap B$.
\item
$\bigcup_{j=1}^{k_{n-1}}B_{n-1,j}=C$.
\end{itemize}

{\it \underline{Step n:}} let $1\leq j\leq k_{n-1}$. Consider
$A_{n,j,i_1},\cdots, A_{n,j,i_{l_{n,j}}}$, the sets in $\mathcal{R}$ such
that
$$
w_n^{-1}.A_{n,j,i_s}\cap B_{n-1,j}\neq \emptyset,  \mbox{ for every }
1\leq s\leq l_{n,j}.$$
We denote by $B_{n,1},\cdots, B_{n,k_n}$, with $k_n=\sum_{j=1}^{k}l_{n,j}$, the elements in the collection
$$\{w_n^{-1}.A_{j,i_s}\cap B_{n-1,j}:1\leq s\leq l_{n,j}, 1\leq j\leq k_{n-1}\,\}.$$ We have
 $\mathcal{Q}_n=\{w.B_{n,j}: w\in
D,\, 1\leq j\leq k_{n}\}$  is a clopen finite
partition of $X$. In addition,  for every  $1\leq l\leq k_n$  there exist $1\leq j\leq k_{n-1}$ and $1\leq s\leq l_{n,j}$
such that $B_{n,l}\subseteq B_{n-1,j}$ and  $B_{n,l}\subseteq w_n^{-1}.A_{n,j,s}$.  This implies that for every $1\leq
i\leq n-1$, $w_i.B_{n,l}\subseteq w_i.B_{n-1,j}$ and  by hypothesis,
$w_i.B_{n,l}$ is included in a subset $A_i$ in $\mathcal{R}$.   Since
$\bigcup_{l=1}^{k_n}B_{n,l}=\bigcup_{i=1}^{k_{n-1}}B_{n-1,i}$, the partition $\mathcal{Q}_n$ satisfies
\begin{itemize}
\item
For every $1\leq j\leq k_n$ and  $1\leq i\leq n$, there exist $A$ in
$\mathcal{R}$ and $B$ in $\mathcal{Q}$ such that $w_i.B_{n,j}$ is
contained in $A\cap B$.
\item
$\bigcup_{j=1}^{k_n}B_{n,j}=C.$
\end{itemize}
This implies that at the end of the step $l$, we obtain a partition
$$\mathcal{R}\wedge\mathcal{Q}=\mathcal{Q}_l=\{w.B_{l,j}: w\in D, \,
1\leq j\leq k_l\},$$ which is finer than $\mathcal{R}$ and
$\mathcal{Q}$, and which satisfies $\bigcup_{j=1}^{k_l}B_{n,j}=C$.

\medskip
Now we use this algorithm to prove the Proposition.
>From Corollary \ref{caracterizacion}, there exists a decreasing
sequence $(C_n)_{n\geq 0}$ of clopen subsets of $X$ and a
decreasing sequence $(\Gamma_n)_{n \geq0}$ of syndetic subgroups of
$G$ such that $|\bigcap_{n\geq 0}C_n|=1$ and $T_{C_n}(x)=\Gamma_n$
for all $x\in C_n$.\\
For every $n\geq 0$, we take a subset $D_n$ of $G$ containing
exactly one representing element in each class of $G/\Gamma_n$, and
we define
  $$\mathcal{Q}_n=\{w.C_n: w\in D_n\}.$$
  The collection  $\mathcal{Q}_n$ is a finite clopen partition of
  $X$.\\
Since $X$ is a Cantor set, it is always possible to take a sequence
$(\mathcal{R}_n)_{n\geq 0}$ of finite clopen partitions of $X$
which spans its topology.\\
We construct the desired  sequence
$(\mathcal{P}_n)_{n\geq 0}$ as follows:
\begin{itemize}
\item
  We set $\mathcal{P}_0=\mathcal{R}_0\wedge\mathcal{Q}_0$.
\item
  For $n>0$. First, we set
  $\mathcal{P}'_n=\mathcal{R}_{n}\wedge\mathcal{Q}_n$, and then
$\mathcal{P}_n=\mathcal{P}_{n-1}\wedge\mathcal{P}_n'$.
\end{itemize}
 From this construction we get
$$
(\mathcal{P}_n=\{w.C_{n,j}: w\in D_n, \, 1\leq j\leq k_n\})_{n\geq
0},$$ a sequence of finite clopen partition of $X$ satisfying, for
every $n\geq 0$:
\begin{enumerate}
\item[(i)]
$\mathcal{P}_n$ is finer than $\mathcal{P}_{n-1}$ and
$\mathcal{R}_n$.
\item[(ii)]
$\bigcup_{j=1}^{k_n}C_{n,j}=C_n$.
\end{enumerate}
The first point implies  $(\mathcal{P}_n)_{n\geq 0}$ is a nested
sequence and that it spans the topology of $X$. The second point
implies that this sequence verifies conditions 1. and  2. from the
Proposition.
\end{proof}

The group $G$ is amenable if and only if it has a F{\o}lner
sequence, that is,  a sequence $(F_n)_{n\geq 0}$ of finite subsets
of $G$ such that for every $g\in G$
$$
\lim_{n\to \infty}\frac{|gF_n\bigtriangleup F_n|}{|F_n|}=0.
$$

\begin{remark}
\label{amenable}
{\rm From \cite{We} we deduce that in the case $G$ amenable, it is
always possible to take $(D_n)_{n\geq 0}$, defined as in Proposition \ref{KRpartition}, as a F{\o}lner sequence.
Thus in this paper, if $G$ is amenable, we suppose that
$(D_n)_{n\geq 0}$ is  F{\o}lner.}
\end{remark}

Let $(X,G)$ be a regularly recurrent Cantor system. Consider the
sequence of finite clopen partitions of $X$ as in Proposition
\ref{KRpartition}:
$$
(\mathcal{P}_n=\{w.C_{n,k}: w\in D_n, \, 1\leq k\leq k_n\})_{n\geq
0}.$$ Let $n\geq 0$. The incidence matrix between $\mathcal{P}_n$
and $\mathcal{P}_{n+1}$ is $A_n\in \mathcal{M}_{k_n\times
k_{n+1}}(\mathbb{Z}^+)$  defined by
$$
A_n(i,j)=|\{w\in D_{n+1}: w.C_{n+1,j}\subseteq C_{n,i}\}|.
$$
Notice that $\sum_{i=1}^{k_n}A_n(i,j)=q_{n,j}$ is the number of $w\in
D_{n+1}$ such that $w.C_{n+1,j}\subseteq C_n$. Since the set of
return times of the points in $C_n$ to $C_n$ is equal to $\Gamma_n$,
the number $q_{n,j}$ does not depend on $j$ and it is equal to the
number of $w\in D_{n+1}$ which are in $\Gamma_n$. So,
$q_{n,j}=\frac{|D_{n+1}|}{|D_n|}$ for every $1\leq j\leq k_{n+1}$.
Consider the set $$\triangle_n=\{(x_1,\cdots, x_{k_n})\in
(\mathbb{R}^+)^{k_n}: \sum_{i=1}^{k_n}x_i=\frac{1}{|D_n|}\}.$$
Since, for every $1\leq j\leq k_{n+1}$,
$\sum_{i=1}^{k_n}A_n(i,j)=\frac{|D_{n+1}|}{|D_n|}$, the map
$A_n:\triangle_{n+1}\to \triangle_n$ is well defined.

\medskip

Because  $(\P_n)_{n\geq 0}$ is a countable collection of clopen
sets that spans the topology of $X$, any invariant measure defined
on this family of sets extends to a unique invariant measure on the
Borel $\sigma$-algebra of $X$. So,
any  invariant measure $\mu$ on $(\mathcal{P}_n)_{n\geq 0}$ must
verify
$$
\mu(C_{n,i})=\sum_{j=1}^{k_{n+1}}A_n(i,j)\mu(C_{n+1,j}), \mbox{ for
every $1\leq i\leq k_n$ and $n\geq 0$ },$$ and it is completely
determined by this relation. In other words, we can identify an
invariant measure with an element in the inverse limit
$\lim_{\leftarrow n}(\triangle_n, A_n)$. In the next Lemma we
provide a sufficient condition to have the  reciprocal.

\begin{lemma}
\label{medidas} If  $G$ is amenable then $\mathcal{M}_{G}(X)$ is
affine-homeomorphic to $\lim_{\leftarrow n}(\triangle_n, A_n)$.
\end{lemma}

\begin{proof}
Since $G$ is amenable, we can suppose that $(D_n)_{n\geq 0}$ is a
F{\o}lner sequence (See Remark \ref{amenable}).\\
Let $((x_{n,1},\cdots, x_{n,k_n}))_{n\geq 0}$ be an element in
$\lim_{\leftarrow n}(\triangle_n, A_n)$. It defines a probability
measure on $X$ by setting
$$
\mu(u.C_{n,i})=x_{n,i}, \mbox{ for every $1\leq i\leq k_n$, $u\in D_n$ and $n\geq
0$}.$$ To show this measure is invariant it is sufficient to show
that for every $n\geq 0$, $1\leq k\leq k_n$ and $v\in G$,
$\mu(v.C_{n,k})=\mu(C_{n,k})=x_{n,k}$.

\medskip

Let $m>n$ and consider the sets
$$J(m,n,k,l)=\{w\in D_m: w.C_{m,l}\subseteq C_{n,k}\}.$$
$$
J_1(m,n,k,l)=\{w\in J(m,n,k,l): vw\in D_m\} \mbox{ and }
J_2(m,n,k,l)=J(m,n,k,l)\setminus J_1(m,n,k,l).$$

We have
$$v.C_{n,k}=\bigcup_{l=1}^{k_m}\bigcup_{w\in J(m,n,k,l)}vw.C_{m,l},$$
and then
$$\mu(v.C_{n,k})=\sum_{l=1}^{k_m}\sum_{w\in
J(m,n,k,l)}\mu(vw.C_{m,l})= \sum_{l=1}^{k_m}\sum_{w\in
J_1(m,n,k,l)}\mu(vw.C_{m,l})+\sum_{l=1}^{k_m}\sum_{w\in
J_2(m,n,k,l)}\mu(vw.C_{m,l}).$$ Since $\mu(u.C_{m,l})=\mu(C_{m,l})$
for $u\in D_m$, we get
$$
\mu(v.C_{n,k})=\sum_{l=1}^{k_m}
|J_1(m,n,k,l)|\mu(C_{m,l})+\sum_{l=1}^{k_m}\sum_{w\in
J_2(m,n,k,l)}\mu(vw.C_{m,l})$$
$$=\mu(C_{n,k})-\sum_{l=1}^{k_m}\sum_{w\in
J_2(m,n,k,l)}\mu(C_{m,l})+ \sum_{l=1}^{k_m}\sum_{w\in
J_2(m,n,k,l)}\mu(vw.C_{m,l}).
$$
Thus we have
$$
|\mu(v.C_{n,k})-\mu(C_{n,k})|\leq \sum_{l=1}^{k_m}\sum_{w\in
J_2(m,n,k,l)}\mu(C_{m,l})+ \sum_{l=1}^{k_m}\sum_{w\in
J_2(m,n,k,l)}\mu(vw.C_{m,l}).
$$
Because  $J_2(m,n,k,l)\subset \{w\in D_m: vw\notin D_m  \}$, we have
$$
|\mu(v.C_{n,k})-\mu(C_{n,k})|  \leq \sum_{\{w\in D_m:\, vw\notin D_m
\}}\sum_{l=1}^{k_m}\mu(C_{m,l})+ \sum_{\{w\in D_m: \, vw\notin D_m
\}}\sum_{l=1}^{k_m}\mu(vw.C_{m,l})
$$
$$
  \,\,\,\,\,\,\,\, =\sum_{\{w\in D_m:\, vw\notin D_m
\}}\mu\left(\bigcup_{l=1}^{k_m}C_{m,l}\right)+ \sum_{\{w\in D_m: \,
vw\notin D_m \}}\mu\left(\bigcup_{l=1}^{k_m}vw.C_{m,l}\right).
$$

Since  $|\{w\in D_m: vw\notin D_m \}|\leq |v.D_m\bigtriangleup D_m|$
and
$\mu\left(\bigcup_{l=1}^{k_m}C_{m,l}\right)=\mu\left(\bigcup_{l=1}^{k_m}vw.C_{m,l}\right)=\frac{1}{|D_m|}$,
we have
$$
|\mu(v.C_{n,k})-\mu(C_{n,k})| \leq \frac{2|v.D_m\bigtriangleup D_m|
}{|D_m|}.
$$
So, because $(D_n)_{n\geq 0}$ is F{\o}lner, we get
$\mu(v.C_{n,k})=\mu(C_{n,k})$.


\end{proof}

\section{Semicocycles}
The notion of a semicocycle has been extensively used in the theory of one-dimensional Toeplitz flows (see
\cite{D}). In this section it is not used but we develop it for actions of a residually finite discrete group
$G$ for further utility.

We fixe a finite family $S$ of generators  of the group $G$, and we suppose furthermore that this family is
symmetric ($S^{-1}=S$). We can then associate to the group $G$ and to the family $S$ a {\it Cayley graph}.
This graph is defined as follow: the vertices are the elements of $G$ and two elements $g_1, g_2$ of the
group $G$ are related by an edge if and only if there exists a element $s$ of $S$ such that $g_2 = s.g_1 $,
where . is the multiplication in the group $G$. This graph is endowed with the natural metric  of the length
path: the distance between two points is the minimal length of paths going from one point to the other, each
edge counting for a longer one. This induces on $G$ a metric $d$, which is invariant by the multiplication to
the right by any element $\gamma$ of $G$.\\
Recall that for a residually finite group $G$ and a decreasing sequence $(\Gamma_i)_{i\geq 0}$ of of syndetic
subgroup of $G$ with $\bigcap_{i\geq 0} \Gamma_i = \{e\}$, the stabilizer of ${\bf e} = (e_i)_{i\geq 0}$ in
the $G$-subodometer $\overleftarrow{G} =\lim_{\leftarrow n} (G/ \Gamma_n, \pi_n)$ is trivial. This defines an
immersion $\tau$ of $G$ into $\overleftarrow{G}$.

\begin{definition}
Let $\overleftarrow{G}=\lim_{\leftarrow n}(G/\Gamma_n ,\pi_n)$ be  a $G$-subodometer with $\bigcap_{i\geq 0}
\Gamma_i = \{e\}$ and let $K$ be a compact metric space. A function $f:G \to K$ is a {\em semicocycle on
$\overleftarrow{G}$} if it is continuous with respect $\Theta_{\overleftarrow{G}}$, where
$\Theta_{\overleftarrow{G}}$ is the topology on $G$ inherited from $\overleftarrow{G}$ (we identify $\tau(G)$
with $G$).
\end{definition}

The functions $f:G\to K$ may be seen as elements of the topological
dynamical system $(K^{G},G)$, where $K^{G}$ is endowed with the metrizable
product topology and  the left-action of $\gamma\in G$ on
$f=(f(g))_{g\in G}\in K^{G}$ is the shift action: this means
$\gamma .f=\{f'\}_{g\in G}$, where $f'(g)= f(g\gamma )$ for every $g \in G$.\\
The proofs of Theorems \ref{semicociclo1} and \ref{semicociclo2}
below follow the same ideas as used in \cite{D} for dimension one.

\begin{theorem}
\label{semicociclo1} If $f\in K^{G}$ is a semicocycle  on some
subodometer $\overleftarrow{G}$ then $f$ is a regularly recurrent
point of $(K^{G},G)$.
\end{theorem}
\begin{proof}
Fix $\epsilon > 0$ and a finite set $C$ in $G$. The pair $(\epsilon, C)$ determines a basic open
set $V$ in the Tychonov topology. Since $f$ is continuous on $G$ for the topology induced by the
odometer $\overleftarrow{G}$, there exists $\delta >0$ such that for every $g \in C$ and $g' \in
G$, ${\rm dist}(g,g') < \delta$ (for the metric inherited from $\overleftarrow{G}$) implies ${\rm
d}( f(g),f(g'))< \epsilon$ in K. By definition of a subodometer, there exist a finite index
subgroup $\Gamma$ of $G$ and a factor map $\pi: \overleftarrow{G} \to G/\Gamma$ such that for any
element $w$ of $G/ \Gamma$, $\pi^{-1}(w)$ is a clopen subset of $\overleftarrow{G}$ with diameter
smaller than $\delta$. Furthermore, for any $y \in \pi^{-1}(w)$, $T_{\pi^{-1}(w)}(y)$ is a group
conjugated to $\Gamma$. Let us consider now the finite index normal subgroup $H= \cap_{g\in G} g
\Gamma g^{-1}$. Since $\Gamma$ is of finite index in $G$, there is just a finite number of groups
conjugated to $\Gamma$ and the former intersection is a finite intersection. The group $H$ is a
subgroup of any group of the kind $T_{\pi^{-1}(w)}(y)$ with $w \in G/ \Gamma, y \in \pi^{-1}(w)$.
Thus, ${\rm dist}(n'.g,n') < \delta$ for any $g \in H$ and hence $d(f(n'.g), f(n'))< \epsilon$ for any $g \in H$ by the normality of $H$. We
prove by this way that the  $H$-orbit of $f$ is included in $V$ and then $f$ is a regularly
recurrent point of $K^G$.
\end{proof}

Proposition \ref{almost} and Theorem \ref{semicociclo1} imply that $(\Omega_{G}(f),G)$ is a minimal
almost 1-1 extension of some subodometer, where $\Omega_{G}(f)$ represents the closure orbit of a
semicocycle $f$ in $K^{G}$ with a trivial stabilizer under the action of $G$. Notice that $\overleftarrow{G}$ need not to be the maximal
equicontinuous factor of $(\Omega_{G}(f),G)$, as we will see later.

\medskip

Let $f\in K^{G}$ be a semicocycle on a $G$-subodometer $\overleftarrow{G}$. Since we have identified the
group $G$ with $G$ embedded in $\overleftarrow{G}$, it makes sense to define $F$ to be the closure of the
graph of $f$ in $\overleftarrow{G}\times K$ endowed with the product topology, $F=\overline{\{(g,f(g)):g\in
G\}}\subseteq \overleftarrow{G}\times K$. Let $F(\textbf{g})$ be the set $\{k\in K: (\textbf{g},k)\in F\}$
for $\textbf{g}\in
\overleftarrow{G}$.\\
We call $C_f$ the set of $\textbf{g}\in \overleftarrow{G}$ such that
$|F(\textbf{g})|=1$ and $D_f=\overleftarrow{G}\setminus C_f$. Since
$f$ is continuous we have that $F(g)=\{g\}$ for all
$g\in G$. Thus $C_f$ is the subset where $f$ can be continuously extended by
$f(\textbf{g})=F(\textbf{g})$.\\
The semicocycle $f$ is said to be {\em invariant under no rotation} if for every ${\bf h}_1 \neq
{\bf h}_2 \in \overleftarrow{G}$ there exists a $g \in G$ such that $F(g.\textbf{h}_1)\neq
F(g.\textbf{h}_2)$.
\begin{theorem}
\label{semicociclo2} Let $(X,G)$ be a minimal topological dynamical system and $\overleftarrow{G}=
\lim_{\leftarrow n}(G/\Gamma_n ,\pi_n)$ be a $G$-subodometer with $\bigcap_{i\geq 0} \Gamma_i = \{e\}$. There
exists an almost 1-1 factor $\pi$ of $(X,G)$ onto $(\overleftarrow{G},G)$ with $|\pi^{-1}({\bf e})| =1$ if
and only if $(X,G)$ is conjugated to $(\Omega_{G}(f),G)$, where $f$ is a semicocycle on $\overleftarrow{G}$,
invariant under no rotation.
\end{theorem}

\begin{proof}
Consider the system $(\Omega_G (f), G)$. By definition, for every $x
\in \Omega_G (f)$, there exists a sequence $(g_i)_i \subset G$ such
that for each $h \in G$, $\lim_i f (h g_i) = x(h)$. Let ${\bf j} \in
\overleftarrow{G}$ be an accumulation point of the sequence
$(g_i)_i$. We have $x(h)\in F(h.{\bf j})$. By a straightforward
calculation, we check that for each such $\bf j$, the set $\{(h.
{\bf j}, x(h)) | \ h\in G\}$ is a dense subset of $F$. Since $f$ is
invariant under no rotation, $\bf j$ is determined for any $x$ in
an unique way.
So we have proved that if $g_i .f \to x$ then $g_i \to {\bf j}$. The map $\pi: x\in \Omega_G(f)
\mapsto {\bf j} \in \overleftarrow{G}$ is a continuous extension onto $\Omega_G(f)$ of the
application $ g.f \mapsto g$. It is straightforward to check that $\pi$ is a factor map that sends
$f$ to $\bf e$. If $\pi(x)= {\bf e}$ then $x(h) \in F(h)=\{f(h)\}$ and $x(h)=f(h)$ for any $h \in
G$. Since the system $(\Omega_G (f), G)$ is minimal, $\pi$ is an almost 1 to 1 factor map.

\noindent Conversely, consider a minimal almost 1-1 extension $(X, G)$ of a $G$-subodometer and
$\pi :X \to \overleftarrow{G}$ the associated factor map. Consider $x \in X$ such that $\pi(x)$ has
a singleton fiber by $\pi$ so this is the same for all the elements of its $G$-orbit. The map $f: g
\in G \mapsto \pi^{-1}(g.{\pi(x)})=g.\pi^{-1}(x) \in X$ is continuous for the induced topology on
$G$, it is then a semicocycle. This is straightforward to check that $F({\bf j})= \pi^{-1}({\bf
k})$ where ${\bf k} \in \overleftarrow{G}$ is the limit point of the sequence $(g_i.\pi(x))_i$ with
$(g_i)$ a sequence of $G$ that converges to ${\bf j}$. The set $\pi^{-1}(\bf k)$ does not depend of
the choice of the sequence $(g_i)$. It is then straightforward to show that $f$ is invariant under
no rotation. The conjugating map from $(\Omega_G(f),G)$ onto $(X,G)$ is the projection onto the
neutral element coordinate: $\phi \mapsto \phi({\bf e})$. By a standard way, we check this
application is a homeomorphism which commutes with the $G$-action.
\end{proof}

\begin{corollary}
A topological dynamical system $(X,G)$ is a minimal almost 1-1 extension of a free odometer
$(\overleftarrow{G},G)$ if and only if it is conjugated to $(\Omega_{G}(f),G)$, where $f$ is a semicocycle on
$G$, invariant under no rotation.
\end{corollary}
\begin{proof}
For a factor map $p: X \to \overleftarrow{G}$ and any point $x\in
X$, by a right multiplication by $\pi(x)^{-1}$, we obtain again a
factor map that sends the point $x$ to ${\bf e}$. The result follows
from Theorem \ref{semicociclo2}.
\end{proof}

\section{$G$-Toeplitz Arrays}

In this section, we suppose that $G$ is a discrete finitely generated group. Let $\Sigma$ be a finite
alphabet and $\Gamma\subseteq G$ a syndetic subgroup
of $G$.  For\\
$x=(x(g))_{g\in G}\in \Sigma^{G}$ we define:
$$
Per(x,\Gamma,\sigma)=\{g\in G: x(g\gamma )=\sigma \, \mbox{ for all }
\gamma\in
\Gamma \}, \,
\sigma \in \Sigma,
$$
$$
Per(x,\Gamma)=\bigcup_{\sigma \in \Sigma}Per(x,\Gamma,\sigma).
$$

\noindent Clearly for two subgroups $\Gamma_1$ and $\Gamma_2$, $\Gamma_1\subset\Gamma_2$, we have
$Per (x,\Gamma_2,  \sigma) \subset
Per(x,\Gamma_1, \sigma)$. When $Per(x,\Gamma)\neq \emptyset$ we say that
$\Gamma$ is a
{\em group of periods of
$x$}. Furthermore, $Per (x,\Gamma)$ is a subset stable by multiplication to
the right
by a element of
$\Gamma$. We say that $x$ is a {\em $G$-Toeplitz array } (or simply a
Toeplitz array)
if for all $g\in
G$ there exists $\Gamma\subseteq G$ syndetic subgroup of $G$ such that $g\in
Per(x,\Gamma)$.
\begin{proposition}
\label{toeplitz} The following statements concerning  $x\in \Sigma^{G}$ are
equivalent:
\begin{enumerate}
\item $x$ is Toeplitz array. \item There exists a sequence of syndetic
subgroups
$(\Gamma_n)_{n\geq0}$,
  such that\\
$\Gamma_{n+1} \subset \Gamma_n$ and $G = \cup_n Per(x,\Gamma_n)$ for
all $n\geq
0$. \item $x$ is
regularly recurrent.
\end{enumerate}
\end{proposition}

\begin{proof}
Let $D_n$ be the ball of radius $n$ in $G$ centered in the neutral element.\\
Suppose that $x$ is a
Toeplitz array. Since for $Z_1$ and $Z_2$, two groups of period of $x$,
we have
$Per(x,Z_1) \subset
Per(x, Z_1\cap Z_2)$, for any $n \geq 0$, there exists  a syndetic
subgroup $Z_n$
such that $D_n
\subset Per(x,Z_n)$. Let $\Gamma_0= Z_0$ and $\Gamma_{n+1}= \Gamma_n \cap
Z_n$. The sequence $(\Gamma_n)_n$ satisfies the statement (2).\\
Let $(\Gamma_n)_n$ be a sequence as in statement (2). Let $C_n$ be the set
$\{y\in \Sigma^{G}:
y(D_n)=x(D_n)\}$ for all $n\geq 0$, $(C_n)_{n\geq 0}$ is a
fundamental system
of clopen
neighborhoods of $x$. Since $D_n$ is contained in   $Per(x,\Gamma_n)$, the set
of return times of
$x$ to $C_n$ contains $\Gamma_n$ which implies that $x$ is regularly
recurrent.\\
Suppose that $x$ is regularly recurrent. For $n\geq 0$ we take $\Gamma_n$ a
syndetic subgroup of
$G$ such that $\Gamma_n\subseteq T_{C_n}(x)$. It holds that $G$ is equal to
$\bigcup_{n\geq
0}Per(x,\Gamma_n)$, which means that $x$ is a Toeplitz array.
\end{proof}

\medskip

A subshift $(X,G)$ is a {\em $G$-Toeplitz system} (or simply a Toeplitz
system)
if there exists a
Toeplitz array $x$ such that $X=\Omega_{G}(x)$. From Theorem \ref{almost} and
Proposition
\ref{toeplitz} we conclude that the family of minimal subshifts which
are almost
1-1 extensions of
subodometers coincides with the family of Toeplitz systems.\\
In order to know the maximal equicontinuous factor of a given
Toeplitz system, we will introduce  the concepts of essential group
of periods and period structure.

\begin{definition}
Let $x\in \Sigma^{G}$. A syndetic group $\Gamma\subset G$  is called an {\em essential group of periods of $x$} if
$Per(x,\Gamma,\sigma)\subseteq Per(x, g^{-1}\Gamma g,\sigma)g^{-1}=
Per(g.x, \Gamma,\sigma)$ for every $\sigma\in \Sigma$  implies that
$g\in \Gamma$.
\end{definition}

\begin{lemma}
\label{essentialgroups} If $\Gamma$ is an essential group of periods
of $x$ then every group of periods $\Gamma'$ satisfying
$Per(x,\Gamma)\subseteq Per(x,\Gamma')$ is contained in $\Gamma$.
\end{lemma}
\begin{proof}
Let $\Gamma$ be an essential group of periods of $x$. Suppose that
$\Gamma'$ is a group of periods such that $Per(x,\Gamma)\subseteq
Per(x,\Gamma')$. For $w\in Per(x,\Gamma,\sigma)$ and $g\in \Gamma'$
we have $w\gamma g\in Per( x,\Gamma', \sigma)$ for every $\gamma\in
\Gamma$. This implies that $x(w\gamma g)=g.x(w\gamma)=\sigma$ for
every $\gamma\in \Gamma$, which means that $w\in
Per(g.x,\Gamma',\sigma)$. Because $\Gamma$ is essential, we conclude
that $g\in \Gamma$ and then $\Gamma'\subseteq \Gamma$.
\end{proof}

\begin{remark}{\rm From Lemma \ref{essentialgroups} we deduce that the
family of the essential groups of periods is contained in the family of the groups generated by
essential periods introduced in \cite{Co} for the case $G=\ZZ^d$.}
\end{remark}

In the following Lemma we show the existence of essential groups of periods.
\begin{lemma}
\label{esencial} Let $x\in \Sigma^{G}$. If $\Gamma\subseteq G$ is a
group of periods of $x$ then there exists $K\subseteq G$  an
essential  group of periods of $x$ such that $Per(x,\Gamma)\subseteq
Per(x,K)$.
\end{lemma}
\begin{proof}
Let $\Gamma\subseteq G$ be a group of periods of $x$ and $\Gamma'$
be a syndetic normal subgroup of $\Gamma$. We call $\hat{\Gamma'}$
the set
$$
\bigcup_{g\in G}\{Hg:\, H \mbox{ syndetic subgroup of $G$ such that
} Per(x,\Gamma',\sigma)\subseteq Per(x,g^{-1}Hg,\sigma)g^{-1},
\forall \sigma\in \Sigma\}.
$$
Let $K$ be the group generated by the elements of $\hat{\Gamma'}$. Let $w \in
Per(x,\Gamma',\sigma)$. For any $\gamma \in \Gamma'$ and any $Hg \in \hat{\Gamma'}$, we have
$w\gamma$ belongs to $Per(x,\Gamma',\sigma) \subseteq Per(x,g^{-1}Hg,\sigma)g^{-1}$. This implies
that  for every $hg \in Hg\in \hat{\Gamma'}$ we have $w\gamma hg \in Per(x,g^{-1}Hg,\sigma)$. Since
$\Gamma'$ is a normal subgroup, we get for any $\gamma\in \Gamma'$ and any $hg \in Hg \in
\hat{\Gamma'}$, $x(whg\gamma)=\sigma$, which means that $whg\in Per(x,\Gamma',\sigma)$. Thus we
obtain that for any $h_1g_1, \ldots, h_ng_n$ with $h_ig_i$ belonging to a set in $\hat{\Gamma'}$,
we have $x(w h_1g_1 \ldots h_ng_n)=\sigma$. In other words, $Per(x,\Gamma',\sigma)$ is contained in
$Per(x,K,\sigma)$. So, we have  $Per(x,\Gamma,\sigma) \subseteq Per(x,\Gamma',\sigma) \subseteq
Per(x,K,\sigma)$. If $g\in G$ is such that $Per(x,K,\sigma)\subseteq
Per(g.x,K,\sigma)=Per(x,g^{-1}Kg,\sigma)g^{-1}, \forall \sigma\in \Sigma$, then $Kg$ belongs to
$\hat{\Gamma'}$, which implies that $g$ is in $K$.
\end{proof}

\begin{corollary}
\label{estructura} Let $x\in \Sigma^{G}$ be a Toeplitz array. There
exists a sequence $(\Gamma_n)_{n \geq0}$ of essential group of
periods of $x$  such that $\Gamma_{n+1}\subseteq \Gamma_n$ and
$\bigcup_{n\geq 0}Per(x,\Gamma_n)=G$.
\end{corollary}

\begin{proof}
>From Proposition \ref{toeplitz} (2) we conclude there exists a decreasing sequence $(\Gamma'_n)_{n\geq 0}$ of
syndetic groups of periods of $x$ such that $\bigcup_{n\geq 0}Per(x,\Gamma'_n)=G$. We set $\Gamma_0$ an
essential group of  periods of $x$ such that $Per(x,\Gamma_0')\subseteq Per(x,\Gamma_0)$.  For $n>0$ we set
$\Gamma_n''=\Gamma'_n\cap \Gamma_{n-1}$ which is a syndetic subgroup of $G$, and since $Per(x,\Gamma_{n-1})$
and $Per(x,\Gamma'_n)$ are contained in $Per(x,\Gamma_n'')$, $\Gamma_n''$ is a group of periods of $x$. Thus,
by Lemma \ref{esencial}, there exists an essential group of periods $\Gamma_n$, such that
$Per(x,\Gamma_{n-1})\subseteq Per(x,\Gamma_n'')\subseteq Per(x,\Gamma_n)$. Since $\Gamma_{n-1}$ is an
essential group of periods, from Lemma \ref{essentialgroups} we get $\Gamma_n\subseteq \Gamma_{n-1}$. Because
$\bigcup_{n\geq 0}Per(x,\Gamma'_n)=G$, we deduce  $\bigcup_{n\geq
0}Per(x,\Gamma_n)=G$. 
\end{proof}

\begin{definition}
A sequence of groups as in Corollary \ref{estructura} is called a {\em period
structure} of $x$.
\end{definition}
In the sequel, we will show that from a period structure $(\Gamma_n)_{n\geq
0}$ of a $G$-Toeplitz
array $x$ it is possible to construct a sequence of nested finite clopen
partitions of
$\Omega_{G}(x)$. From this sequence of partitions it will be easy to define an
almost 1-1 factor
map between the Toeplitz system $(\Omega_{G}(x),G)$ and the
odometer $\overleftarrow{G}=\lim_{\leftarrow n}(G/\Gamma_n,\pi_n)$.\\
Let $x\in \Sigma^{G}$ be a Toeplitz array, let $y\in \Omega_{G}(x)$ and let
$\Gamma\subseteq G$ be
a subgroup of $G$ with finite index. Since $(\Omega_{\Gamma}(y),\Gamma)$ is
minimal, if $\Gamma$ is
a group of periods of $y$ then $\Omega_{\Gamma}(y)\subseteq C_{\Gamma}(y)$,
where
$$
C_{\Gamma}(y)=\{x'\in \Omega_{G}(x):
Per(x',\Gamma,\sigma)=Per(y,\Gamma,\sigma),
\, \, \, \forall
\, \, \sigma\in \Sigma\}.
$$

\begin{lemma}
\label{igualdad} $C_{\Gamma}(y)=\gamma.C_{\Gamma}(y)$ for every
$\gamma\in \Gamma$.
\end{lemma}
\begin{proof}
Let $x'\in \gamma.C_{\Gamma}(y)$. There exists $x''\in
C_{\Gamma}(y)$ such that $x'=\gamma.x''$. If $g\in
Per(x',\Gamma,\sigma)$ then $x'(g\gamma')=\sigma$ for every
$\gamma'\in \Sigma$. In particular, we have
$$x'(g\gamma'\gamma^{-1})=\gamma^{-1}.x'(g\gamma')=x''(g\gamma')=\sigma,
\, \, \forall \gamma'\in \Gamma,$$ which implies
$Per(x',\Gamma,\sigma)\subseteq Per(y,\Gamma,\sigma)$. On the other
hand, if $g\in Per(x'',\Gamma,\sigma)$ then
$$x''(g\gamma')=x''(g\gamma'\gamma)=\gamma.x''(g\gamma')=x'(g\gamma')=\sigma,
\,
\,
\forall \gamma'\in \Gamma,$$ which implies that
$Per(y,\Gamma,\sigma)\subseteq Per(x',\Gamma,\sigma)$. Thus we
obtain that $\gamma.C_{\Gamma}(y)\subseteq C_{\Gamma}(y)$. Since
this is true for every $\gamma\in \Gamma$, we conclude that
$\gamma.C_{\Gamma}(y)=C_{\Gamma}(y)$.
\end{proof}

We will use the following convention: For a $\Gamma$-periodic subset
$C$ of $\Omega_{G}(x)$, i.e., such that $w.C = w'.C$ whenever
$w^{-1}w' \in  \Gamma$ we will write $v.C$ instead of $w.C$, where
$v$ is the projection of $w$ to $G/\Gamma$.

\begin{proposition}
\label{particion} Let $x\in \Sigma^{G}$ be a Toeplitz array and let
$y\in \Omega_{G}(x)$. If $\Gamma\subseteq G$ is a subgroup generated by
essential periods of  $y$  then
$\Omega_{\Gamma}(y)=C_{\Gamma}(y)$ and $\{w.C_{\Gamma}(y)\}_{w\in
G/\Gamma}$ is a clopen partition of $\Omega_{G}(x)$.
\end{proposition}
\begin{proof}
By Lemma \ref{igualdad}, $C_{\Gamma}(y)$ is a clopen set and we have
$\Gamma\subseteq T_{C_{\Gamma}}(x')$ for every $x'\in
C_{\Gamma}(y)$. In the sequel, we will show that for $\Gamma$ a
group generated by essential periods, we have
$T_{C_{\Gamma}(y)}(x')=\Gamma$ for every $x'\in C_{\Gamma}(y)$,
which will allow us to conclude.

\medskip

Suppose that $g\in \Gamma$ satisfies $g.y\in C_{\Gamma}(y)$. This
implies $Per(g.y,\Gamma,\sigma)=Per(y,\Gamma,\sigma)$ for every
$\sigma\in \Sigma$. Since $Per(g.y,\Gamma,\sigma)=Per(y,g^{-1}\Gamma
g,\sigma)g^{-1}$, we obtain $g\in \Gamma$ because $\Gamma$ is a
group generated by essential periods of $y$. By minimality, we
conclude that $ T_{C_{\Gamma}(y)}(x')=\Gamma$ for every $x'\in
C_{\Gamma}(y)$.  Thus we get that $\{w.C_{\Gamma}(y)\}_{w\in
G/\Gamma}$ is a collection of disjoint sets. Moreover, this
collection is a partition of $\Omega_{G}(x)$ because
$w.\Omega_{\Gamma}(y)\subseteq w.C_{\Gamma}(y)$ for every $w\in
G/\Gamma$ and $\{w.\Omega_{\Gamma}(y)\}_{w\in G/\Gamma}$ is a
covering of $\Omega_{G}(x)$. This also implies that
$\Omega_{G}(x)=C_{\Gamma}(x)$.
\end{proof}

\begin{proposition}
\label{eqfactor} Let $x\in \Sigma^{G}$ be a Toeplitz array. If
$(\Gamma_n)_{n \geq0}$ is a period
structure of $x$ then the subodometer $\overleftarrow{G}=\lim_{\leftarrow
n}(G/\Gamma_n,\pi_n)$ is the maximal
equicontinuous factor of $(\Omega_{G}(x),G)$.
\end{proposition}
\begin{proof}
By Proposition \ref{particion}, if $(\Gamma_n)_{\geq 0}$ is period structure
of the Toeplitz array
$x$, then $(C_{g.\Gamma_n}(x): g\in G/\Gamma_n)_{n\geq 0}$ is a sequence of
nested clopen partitions of
$\Omega_{G}(x)$. This implies that the function $f_n:\Omega_{G}(x)\to
G/\Gamma_n$ given by
$f_n(y)=g$ if and only if $y\in g.C_{\Gamma_n}(x)$ is a well defined
continuous
function, $y\in
\Omega_{G}(x)$, $n\geq 0$. The function $\pi:\Omega_{G}(x)\to
\overleftarrow{G}$
given by $\pi=(f_n)_{n\geq
0}$ is a factor map. Since $\bigcap_{n\geq 0}C_{\Gamma_n(x)}=\{x\}$, we have
that
$\pi^{-1}\{\textbf{e}\}=\{x\}$ and then $\pi$ is an almost 1-1 factor map.
\end{proof}

\begin{theorem}\label{01}
For every subodometer  $\overleftarrow{G}$  there exists a Toeplitz
array $x\in \{0,1\}^{G}$ such that $\overleftarrow{G}$ is the
maximal equicontinuous factor of $(\Omega_{G}(x),G)$.
\end{theorem}
\begin{proof}
Let $\overleftarrow{G}=\lim_{\leftarrow n}(G/\Gamma_n,\pi_n)$ be a
subodometer with $\Gamma_0 =
G$. We distinguish two cases:\\
{\em Case 1:} There exists $m\geq 0$ such that $\Gamma_n=\Gamma_m$
for all $n\geq m$. In this case $\overleftarrow{G}$ is the finite
group $G/\Gamma_m$ and then every minimal almost 1-1 extension will
be conjugate to $\overleftarrow{G}$. For example, $x\in \{0,1\}^{G}$
defined by $x(v)=0$ for all $v\in \Gamma_m$ and $x(v)=1$ if not,
provides a Toeplitz sequence $x$ such that $\overleftarrow{G}$ is
the maximal equicontinuous factor of the system associated to $x$.

\medskip

\noindent {\em Case 2:} For every $m\geq 0$ there exists $n>m$ such that $\Gamma_n\neq \Gamma_m$. In this
case we can take a subsequence $(\Gamma_n)_{n \geq0}$ such that $\Gamma_{n+1}\neq \Gamma_n$ and $[\Gamma_{n}:
\Gamma_{n+1}]\geq 2$ for all $n\geq 0$. By Proposition \ref{factor}, $\overleftarrow{G}$ is conjugate to the
subodometer obtained from this sequence. In order to construct the Toeplitz array $x$ we will consider a
sequence $(D_n)_{ n\geq 0}$  of compact subsets of $G$ such that:
\begin{itemize}
\item for each $n$, $D_n$ is a fundamental domain of $\Gamma_n$ (i.e. $D_n$ contains an unique
element of each class of $G/\Gamma_n$). The set $D_0$ is the singelton set $\{e\}$.

\item For each $n$, $D_n \subset D_{n+1}$ and $D_{n+1}=\bigsqcup_{k\in K_n} D_n.k$ for some finite
set $K_n \subset G$ containing the neutral element $e$ of $G$. By assumption, the cardinal of $K_n$
  is bigger than $2$.

\item $\bigcup_{n \geq 0} D_n =G$.
\end{itemize}
\noindent We define now a sequence of subsets of $G $  $(S_n)_{n
\geq 0}$  by induction. Let $S_0$ be the singleton  $\{e\}$. Let
$v_{1}$ be an element of $D_{1}$ distinct from $e$ and let
$S_{1}=\{v_{1}\}$. For $n>1$, let  $S_n$ be the set $v_{n-1}.
\Gamma_{n-1}\cap D_{n}\setminus D_{n-1}$ and let $v_{n}$ be a point
in $S_{n}$. We define then $x\in \{0,1\}^{G}$ by :
\begin{equation}
  x(w)=\left\{ \begin{array}{ll}
       0 &  \mbox{if $w$ belongs  to $\displaystyle{\cup_{n\geq0}
S_{2n}.\Gamma_{2n+1} } $}\\
                   1 &  \mbox{else}
                 \end{array}
                 \right.
\end{equation}
Remark that $x(w)=1$ for the element $w$ of
$\displaystyle{\cup_{n\geq0} S_{2n+1}.\Gamma_{2n+2} }$. Since
$\bigcup_{j=0}^{n-1} v_j.\Gamma_{j+1} \subseteq Per(x,\Gamma_n,0)$
and $(D_{n-1} \setminus\bigcup_{j=0}^{n-1} v_j.\Gamma_{j+1})
\subseteq Per(x,\Gamma_n,1)$, it holds that $G=\bigcup_{n\geq
0}Per(x,\Gamma_n)$ and $x$ is a Toeplitz array. To conclude that
$\overleftarrow{G}$ is the maximal equicontinuous factor of the
system associated to $x$, by Proposition \ref{eqfactor}, it suffices
to show that $(\Gamma_n)_{n \geq0}$ is a period
structure of $x$.\\

Let us prove by induction on $n$ that $\Gamma_{n}$ is a group generated by essential periods of
$x$. For $n=0$, $\Gamma_{0}=G$ and this is obviously true. Suppose now that $n>0$ and that
$\Gamma_{n-1}$ is a group generated by essential periods. Let $g\in G$ be such that
$Per(x,\Gamma_{n},\sigma) \subset Per(x,g^{-1}\Gamma_{n}g,\sigma).g^{-1}$, for all $\sigma$ of
$\{0,1\}$. Since $\Gamma_{n}\subset \Gamma_{n-1}$, we have $Per(x,\Gamma_{n-1},\sigma) \subset
Per(x,\Gamma_{n},\sigma)$. Let $w$ be in $Per (x, \Gamma_{n-1}, \sigma)$ and $\gamma_{n-1}$ in
$\Gamma_{n-1}$, there exists $\gamma \in D_{n} $ and $\gamma_{n} \in \Gamma_{n}$ such that
$\gamma_{n-1}= \gamma \gamma_{n}$. Then $w\gamma_{n-1}\gamma_{n}^{-1}=w\gamma$ belongs to $Per
(x,\Gamma_{n-1},\sigma) \subset Per(x,g^{-1}\Gamma_{n}g,\sigma).g^{-1} = Per(g.x,
\Gamma_{n},\sigma)$. So we have $\sigma=g.x(w\gamma)=g.x(w \gamma. \gamma_{n})=
g.x(w.\gamma_{n-1})$ and therefor $w \in Per(g.x, \Gamma_{n-1},\sigma) =
Per(x,g^{-1}\Gamma_{n-1}g,\sigma).g^{-1}$ for all $w \in Per (x, \Gamma_{n-1}, \sigma)$. By the
hypothesis of induction we get that $g$ belongs to
$\Gamma_{n-1}$.\\
By the definition of $x$, the element $v_{n-1}$ belongs to $Per (x,
\Gamma_{n}, \sigma)$ with $\sigma = x(v_{n-1})$, so $x(v_{n-1}.g) =
\sigma$. Since $g\in \Gamma_{n-1}$ and by the construction of $x$,
$g$ belongs to $\Gamma_{n}$ and so  $\Gamma_{n}$ is a group
generated by essential periods of $x$.
\end{proof}

\begin{remark}
It is interesting to note that when $\overleftarrow{G}$ is a free odometer, the action of $G$ on
$\overleftarrow{G}$ is free and minimal. The $G$-Toeplitz array $x$, constructed as above, is such that
$(\Omega_G(x), G)$ is an extension almost 1-1 of the system $(\overleftarrow{G}, G)$, so the action of $G$ on
$\Omega_G (x)$ is also free and minimal. All the elements of $\Omega_G(x)$ are not stable by the action of
$G$. Remark also that, very recently and independently from our work, F. Krieger in \cite{Kr} gives a similar
construction of sequence $G$-Toeplitz. This kind of examples are, at our knowledge, the only examples given
with these properties for a general $G$-action.
\end{remark}

\end{document}